\input amstex
\documentstyle{amsppt}
\nologo
\NoBlackBoxes  

\title Equations Defining Toric Varieties \endtitle
\author Bernd Sturmfels\endauthor
\leftheadtext{BERND STURMFELS}

\address Department of Mathematics, University of California,
Berkeley, CA 94720 \endaddress
\email bernd\@math.berkeley.edu\endemail

\subjclass Primary 14M25, 52B20; Secondary 13P10
\endsubjclass

\thanks Supported in part by the NSF and the
David and Lucile Packard Foundation.\endthanks

\endtopmatter

\document

This article is based on my lecture at the
AMS Summer Institute in Algebraic Geometry at Santa Cruz,
July 1995. The topic is toric ideals, by which I mean
the defining ideals of subvarieties of affine or projective space which 
are parametrized by monomials. With one exception
there will be no proofs given in this paper. For most assertions 
which are stated without a reference, the reader can find a proof and 
combinatorial details in my monograph ``Gr\"obner Bases and
Convex Polytopes'' \cite{Stu}. Connections to recent developments
in algebraic geometry and theoretical physics are discussed
in the beautiful survey by David Cox \cite{Cox}.

\vskip 1cm

\head 1. Toric ideals and their polytopes \endhead

In many branches of mathematics and its applications
one encounters algebraic varieties which are parametrized 
by monomials. Such varieties are called {\it toric varieties}
in this article. This stands in contrast to common practise in algebraic
geometry (see \cite{Cox}), where toric varieties are assumed to be
normal. From the point of view  taken by the author
it is more natural to  start out with the following
definition.
Let ${\Cal A} = ( {\bold a}_1,\ldots,{\bold a}_n )$
be any integer $d \times n$-matrix.
Each column vector ${\bold a}_i = (a_{1i},\ldots,a_{di})^T$
is identified with a Laurent monomial
$\,{\bold t}^{{\bold a}_i} =  t_1^{a_{1i}} \cdots t_d^{a_{di}}$.
The {\it toric ideal} $I_{\Cal A}$ associated with ${\Cal A}$
is the kernel of the $k$-algebra homomorphism
$$  {\Bbb C}[x_1,x_2,\ldots,x_n] \rightarrow
{\Bbb C}[t_1,\ldots,t_d,t_1^{-1},\ldots,t_d^{-1}],
\, x_i \,\mapsto \,{\bold t}^{{\bold a}_i} .$$
Every vector ${\bold u}$ in ${\Bbb Z}^n$ can be written uniquely as
${\bold u} = {\bold u}_+ -  {\bold u}_-$, where
${\bold u}_+ $ and $ {\bold u}_-$ are non-negative and have
disjoint support. 
The difference of two monomials is called a {\it binomial}.
An ideal generated by binomials
is called a {\it binomial ideal}.

\proclaim{Lemma 1.1} 
\roster
\item"(a)"  The toric ideal $I_{\Cal A}$
is generated by the binomials
${\bold x}^{{\bold u}_+} - 
 {\bold x}^{{\bold u}_-}$,
where ${\bold u}$ runs over all integer vectors
in the kernel of the matrix ${\Cal A}$.
\item"(b)"  An ideal in $ {\Bbb C}[x_1,\ldots,x_n]$
is toric if and only if it is prime and binomial.
\endroster
\endproclaim

Here is an easy method for computing generators of $I_{\Cal A}$. 
Assume  for simplicity that  ${\Cal A} \subset {\Bbb N}^d$.
Then our toric ideal equals the elimination ideal
$$ I_{\Cal A} \quad = \quad
\langle \,  x_1 - {\bold t}^{{\bold a}_1}, \,
 x_2 - {\bold t}^{{\bold a}_2}, \, \cdots \,
 x_n - {\bold t}^{{\bold a}_n} \,\rangle
\,\,\cap \,\,  {\Bbb C}[x_1,x_2,\ldots,x_n],  \tag 1.1 $$
which can be computed by lexicographic Gr\"obner bases
in ${\Bbb C}[x_1,\ldots,x_n,t_1,\ldots,t_d]$.
More efficient algorithms for the same task are
described in Section 12.1 of \cite{Stu}.

The zero set of $I_{\Cal A}$ in affine $n$-space
is denoted $X_{\Cal A}$ and called the 
{\it affine toric variety} defined by ${\Cal A}$.
The dimension of $X_{\Cal A}$  equals the rank of 
the matrix ${\Cal A}$. If all columns of ${\Cal A}$ have
the same coordinate sum, then the ideal $I_{\Cal A}$ is
homogeneous and defines a 
{\it projective toric variety} $Y_{\Cal A}$ in ${\Bbb P}^{n-1}$.
In what follows we identify ${\Cal A}$ with the
point configuration given by its columns.

\example{Examples 1.2}
Here are some familiar examples of projective toric varieties.
\roster
\item"(a)" The {\it twisted cubic curve} in ${\Bbb P}^3$ is defined
by four equidistant points on a line:
$$ {\Cal A} \quad = \quad
\pmatrix
3 & 2 & 1 & 0 \\
0 & 1 & 2 & 3 \\
\endpmatrix $$
The corresponding toric ideal is generated by three
quadratic binomials:
$$ I_{\Cal A} \quad = \quad
\langle   x_1 x_3 - x_2^2, \,
 x_1 x_4 - x_2 x_3, \,  x_2 x_4 - x_3^2 \rangle. $$
More generally, the {\it
$r$-th Veronese embedding of ${\Bbb P}^{n-1}$} 
equals $Y_{\Cal A}$ for
$$ {\Cal A} \quad = \quad
\bigl\{ \,(i_1,i_2,\ldots,i_n) \in {\Bbb N}^n \,:\,
i_1 + i_2 + \cdots + i_n = r \, \bigr\}.$$
\item"(b)" All {\it rational normal scrolls} are toric.
For instance, the {\it cubic scroll} $S_{2,1}$ in ${\Bbb P}^4$
is defined by the matrix
$$ {\Cal A} \quad = \quad
\pmatrix
2 & 1 & 1 & 0 & 0 \\
0 & 1 & 0 & 2 & 1 \\
0 & 0 & 1 & 0 & 1 \\
\endpmatrix $$
If we were to add a column vector $(0,0,2)$ to 
this matrix then we would get the
quadratic Veronese embedding of ${\Bbb P}^2$ into ${\Bbb P}^5$.
\item"(c)" The Segre embedding of
${\Bbb P}^{r} \times {\Bbb P}^s$ into
${\Bbb P}^{rs+r+s}$ is toric.
Here $I_{\Cal A}$ is the ideal of $2 \times 2$-minors of
an $(r+1) \times (s+1)$-matrix of indeterminates, and
the configuration ${\Cal A}$ consists
of the $rs$ vertices of the product of two regular simplices
$\Delta_r \times \Delta_s$. 
\item"(d)" Consider a generic point  in the
Grassmannian of lines in ${\Bbb P}^3$.
The closure of its orbit under the
natural $({\Bbb C}^*)^4$-action is the toric variety
$Y_{\Cal A}$, where
$$ {\Cal A} \quad = \quad
\pmatrix
1 & 1 & 1 & 0 & 0 & 0 \\
1 & 0 & 0 & 1 & 1 & 0 \\
0 & 1 & 0 & 1 & 0 & 1 \\
0 & 0 & 1 & 0 & 1 & 1 \\
\endpmatrix $$
The closure of any torus orbit in a
flag variety arises from a configuration ${\Cal A}$ of weights 
in a $GL_n({\Bbb C})$-module.
That module is $\wedge^2 {\Bbb C}^4$ in the above case. \qed
\endroster
\endexample

\vskip .1cm
\noindent
The $d$-dimensional algebraic torus $({\Bbb C}^*)^d $ acts 
on affine $n$-space ${\Bbb C}^n$ via
$$ (x_1,\ldots,x_n) \,\,\, \mapsto \,\,\,
(x_1 {\bold t}^{{\bold a}_1},\ldots,x_n {\bold t}^{{\bold a}_n}).$$
The affine toric variety $X_{\Cal A}$ is the closure of the
$({\Bbb C}^*)^d $-orbit of the point $(1,1,\ldots,1)$.
A basic invariant of $X_{\Cal A}$ is the {\it convex polyhedral cone} 
$pos({\Cal A})$
consisting of all non-negative linear combinations of
column vectors in ${\Cal A}$. For a projective
toric variety $Y_{\Cal A}$ we also consider the convex hull
$conv({\Cal A})$ of the points in ${\Cal A}$.
This is a {\it convex polytope} of dimension $\, rank({\Cal A})-1$.
Note that $pos({\Cal A})$ equals the cone over $conv({\Cal A})$.

\proclaim{Proposition 1.3} 
\roster
\item"(a)" The $ ({\Bbb C}^*)^d $-orbits on the affine toric variety
$X_{\Cal A}$ are in order-preserving \break bijection with the faces
of the cone $pos({\Cal A})$.
\item"(b)"  The $ ({\Bbb C}^*)^d $-orbits on the projective toric variety
$Y_{\Cal A}$ are in order-preserving bijection with the faces
of the polytope $conv({\Cal A})$.
\endroster
\endproclaim

We next determine the degree of a projective toric variety $Y_{\Cal A}$
in ${\Bbb P}^{n-1}$.
Let ${\Cal L}$ be the sublattice of ${\Bbb Z}^d$ affinely generated 
by ${\Cal A}$. We normalize the volume form on
${\Cal L} \otimes_{\Bbb Z} {\Bbb R}$ in such a way that
each primitive lattice simplex has unit volume.
The normalized volume of the polytope $conv({\Cal A})$
is a positive integer denoted $\hbox{\it Vol}({\Cal A})$. 

\proclaim{Theorem 1.4}
The degree of $Y_{\Cal A}$ in ${\Bbb P}^{n-1}$
equals the normalized volume $ \hbox{\it Vol}({\Cal A})$.
\endproclaim

It is instructive to verify Theorem 1.4
for the varieties in Examples 1.2.
\roster
\item"(a)" For the twisted cubic, $conv({\Cal A})$
is a line segment of normalized length $3$.
\item"(b)" For the cubic scroll $S_{2,1}$, $conv({\Cal A})$
is a quadrangle with three edges of unit length and
one edge of length two. Its normalized area equals
$\hbox{\it Vol}({\Cal A})=3$.
\item"(c)" The product of simplices $\Delta_r \times \Delta_s$
can be triangulated into ${r+s \choose r}$ simplices of unit volume.
Hence the degree of the Segre variety equals ${r+s \choose r}$.
\item"(d)" For the generic torus orbit on the
Grassmannian of lines in ${\Bbb P}^3$, the polytope
$conv({\Cal A})$ is a regular octahedron. It can
be triangulated into four tetrahedra of unit volume.
Therefore $Y_{\Cal A}$ is a toric threefold of
degree four in ${\Bbb P}^5$.
\endroster

\vskip .1cm

The toric ideals $I_{\Cal A}$ in Example 1.2
are generated by quadratic binomials, and their varieties 
$Y_{\Cal A}$ are projectively normal (see Section 2).
The quadratic generators are easy to find, in view 
of the special structure of the matrices ${\Cal A}$.
The simplicity of these geometric examples is misleading:	
for a general matrix ${\Cal A}$ it is difficult to identify 
generators for $I_{\Cal A}$. One objective of this article is to 
describe methods for studying and solving this problem.
We illustrate this issue for two families of
toric ideals which arise from an application
in computational statistics.

\example{Example 1.5} (The toric variety of the Birkhoff polytope) 
\hfill \break
Let ${\Cal A}$ be the set of $p \times p$-permutation matrices.
Here $\,d = p^2 $, $n = p !$, and 
$conv({\Cal A})$ is the {\it Birkhoff polytope}
consisting of all doubly-stochastic matrices. The dimension 
of $conv({\Cal A})$  and hence of the toric variety $Y_{\Cal A}$
is $(p-1)^2 $. The variables in the
toric ideal $I_{\Cal A}$ are indexed by permutations.
For instance, for $p=3$ we have
$$ I_{\Cal A} \quad = \quad 
\langle \, x_{123} x_{231} x_{312} \,-\,
x_{132} x_{213} x_{321} \,\rangle. $$
In general, the toric ideal $I_{\Cal A}$ is
generated by forms of degree $p$. It is a
challenging combinatorial problem to determine
$\,\hbox{\it Vol}({\Cal A}) = degree(Y_{\Cal A})$.
The known values are:
$$ \matrix
                 p &  &  3 & 4  &  5 & 6 & 7 \\
dim(Y_{\Cal A})    &  & 4 & 9 & 16 & 25 & 36 \\
degree(Y_{\Cal A}) &  &  3 & 352 & 4718075 & 14666561365176 
& 17832560768358341943028 \\
 \endmatrix $$
Another open problem is to bound the
degree of the universal Gr\"obner basis of $I_{\Cal A}$.
\endexample

\example{Example 1.6} (Looks like Segre but isn't)
\hfill \break
Fix integers $r \leq s \leq t$.
Let $n = rst$, $d = rs + rt + st$ and identify $\,{\Bbb Z}^d \,$ with
the direct sum $\,{\Bbb Z}^{r \times 
s} \oplus {\Bbb Z}^{r \times t} \oplus {\Bbb Z}^{s \times t}$.
We denote the standard basis vectors in the three components
as ${\bold e}_{ij}$, ${\bold e}'_{ik}$ and ${\bold e}''_{jk}$ respectively.
Our configuration in this example is
$$ {\Cal A} \quad = \quad \bigl\{\,
 {\bold e}_{ij} \oplus {\bold e}'_{ik} \oplus {\bold e}''_{jk} \,\,:\,\,
 i = 1,\ldots,r, \,\,
 j = 1,\ldots,s, \,\,
 k = 1,\ldots,t \, \bigr\}. $$
The toric ideal $I_{\Cal A}$ is the kernel of the ring map
$$ {\Bbb C}\bigl[ x_{ijk}  \bigr ]\,\rightarrow \,
{\Bbb C} \bigl[ u_{ij}, v_{ik}, w_{jk}  \bigr ] \, , \,\,\,
 x_{ijk} \,\mapsto \,  u_{ij} \cdot v_{ik} \cdot w_{jk}, $$
where the indices $i,j,k$ run as in ${\Cal A}$.
The smallest example is $r=s=t=2$, where
$$ I_{\Cal A} \quad = \quad
 \langle \, x_{111} x_{122} x_{212} x_{221} \, - \,
            x_{112} x_{121} x_{211} x_{222} \,\rangle. $$
At first glance the projective toric variety $Y_{\Cal A}$
looks similar to a Segre variety. But this is a deception.
The following questions are widely open
for general $r,s,t \geq 3$:
\roster
\item"$\bullet$" Characterize the faces of $conv({\Cal A})$,
i.e.~the torus orbits on $Y_{\Cal A}$.
\item"$\bullet$" Determine the normalized volume
of $conv({\Cal A})$,  i.e.~the degree of  $Y_{\Cal A}$.
\item"$\bullet$" Find minimal generators for $I_{\Cal A}$,
or at least bound their degree.
\endroster
\endexample

\vskip 1cm

\head 2. Normal toric varieties \endhead

In algebraic geometry one often assumes that
$X_{\Cal A}$ is normal and that $Y_{\Cal A}$
is projectively normal. This imposes strong combinatorial restrictions
on the configuration ${\Cal A}$. We first discuss these 
restrictions and then we present known results
and open problems concerning the generators of 
the corresponding toric ideals $I_{\Cal A}$.

Let ${\Bbb N}{\Cal A}$ denote the
semigroup spanned by ${\Cal A}$.
Throughout this paper we assume that 
$\, {\Bbb N}{\Cal A} \,\cap \,
 - {\Bbb N}{\Cal A} \,\, =\,\, \{0\} $.
This condition means that the
semigroup algebra $\,{\Bbb C}[{\Cal A}] := 
{\Bbb C}[x_1,\ldots,x_n]/I_{\Cal A}\,$ has no 
non-trivial units. The semigroup ${\Bbb N}{\Cal A}$ lies in 
the intersection of the abelian group ${\Bbb Z}{\Cal A}$ and the 
convex polyhedral cone  $pos({\Cal A})$:
$$ {\Bbb N}{\Cal A} \quad \subseteq \quad
pos({\Cal A}) \,\cap \,{\Bbb Z}{\Cal A}. \tag 2.1 $$

We say that the configuration ${\Cal A}$  
is {\it normal} if equality holds in (2.1).

\proclaim{Lemma 2.1}
The affine variety $X_{\Cal A}$ is normal if and 
only if ${\Cal A}$ is normal.
\endproclaim

If ${\Cal A}$ is not normal, then one can replace it
(e.g.~using Algorithm 13.2 in \cite{Stu})
by the unique minimal finite subset 
${\Cal B}$ of ${\Bbb Z}^d$ such that
$$ {\Bbb N}{\Cal B} \quad = \quad
pos({\Cal A}) \,\cap \,{\Bbb Z}{\Cal A} .$$
In this case $X_{\Cal B}$ is 
the {\it normalization} of $X_{\Cal A}$.

Singularities of toric varieties
are characterized as follows:

\proclaim{Lemma 2.2}
The affine toric variety $X_{\Cal A}$ is smooth if and 
only if the semigroup ${\Bbb N}{\Cal A}$ is isomorphic
to the free semigroup ${\Bbb N}^r$ for some $r$.
\endproclaim

\example{Example 2.3}
The affine toric surface  $X_{\Cal A}$ defined by
${\Cal A} = \{ (2,0),(1,1),(0,2) \} $ is normal but not smooth. 
It is the cone over a smooth quadric curve in ${\Bbb P}^2$.
\endexample

From now on we assume that the ideal $I_{\Cal A}$ is homogeneous.
It defines a projective toric variety $Y_{\Cal A}$ in ${\Bbb P}^{n-1}$.
The intersection of $Y_{\Cal A}$ with the affine chart
$ \{ x_i \not= 0 \} \simeq {\Bbb A}^{n-1}$ equals
the affine toric variety  $ X_{{\Cal A} - {\bold a}_i}$ defined by
the configuration
$$ {\Cal A} - {\bold a}_i \quad = \quad
\bigl\{  {\bold a}_1 - {\bold a}_i, \ldots,
{\bold a}_{i-1} - {\bold a}_i,\,
{\bold a}_{i+1} - {\bold a}_i, \ldots,
{\bold a}_n - {\bold a}_i  \bigr\}.$$
Thus $Y_{\Cal A}$ has an open cover consisting
of $n$ affine toric varieties. In general the 
number $n$ can be lowered, by the following proposition.

\proclaim{Proposition 2.4}
The projective toric variety $Y_{\Cal A}$
is covered irredundently by the affine 
varieties $ X_{{\Cal A} - {\bold a}_i}$ where
${\bold a}_i$ runs over the vertices of the polytope $conv({\Cal A})$.
\endproclaim

Thus $Y_{\Cal A}$ is normal (resp.~smooth) if and only if
each affine chart $ X_{{\Cal A} - {\bold a}_i}$ in the
above cover is normal (resp.~smooth). In particular,
$Y_{\Cal A}$ is normal if and only if
$$ {\Bbb N}( {\Cal A} - {\bold a}_i) \quad = \quad
pos( {\Cal A} - {\bold a}_i) \, \cap \,
{\Bbb Z}( {\Cal A} - {\bold a}_i) \qquad 
\hbox{for \ \ } i=1,\ldots,n .$$
Let $H_{\Cal A}(s)$ denote the {\it Hilbert polynomial} of 
$Y_{\Cal A}$, that is, $H_{\Cal A}(s) = dim_k ({\Bbb C}[{\Cal A}]_s)$
for $s \gg 0$. The {\it Ehrhart polynomial} $E_{\Cal A}$
of the lattice polytope $conv({\Cal A})$ is defined by
$$ E_{\Cal A}(s) \quad= \quad
\# \bigl(\, s \cdot conv({\Cal A}) \,\,\cap \,\,{\Bbb Z}^d \, \bigr)
\qquad \hbox{for all} \,\, s \in {\Bbb N} . $$
The relationship between the Hilbert polynomial and the 
Ehrhart polynomial was studied by A.~Khovanskii in \cite{Kho}.

\proclaim{Proposition 2.5}
A projective toric variety $Y_{\Cal A}$ is normal if and
only if its Hilbert polynomial $H_{\Cal A}$ is equal to
its Ehrhart polynomial $E_{\Cal A}$.
\endproclaim

A much stronger requirement is to ask that
$Y_{\Cal A}$ be {\it projectively normal},
which means that the affine cone $X_{\Cal A}$ 
over $Y_{\Cal A}$ is normal, i.e.,
$\, {\Bbb N}{\Cal A} \,=\,
pos({\Cal A}) \,\cap \,{\Bbb Z}{\Cal A} $.

\example{Example 2.6}
(Normal versus projectively normal)   \hfill \break
Let $d=2,n=4, r \geq 4$ and
$\, {\Cal A}= \{ (r,0), (r-1,1), (1,r-1), (0,r) \}$.
Here $Y_{\Cal A}$ equals the projective line ${\Bbb P}^1$.
It is smooth (hence normal) but not  projectively normal.
The toric ideal $I_{\Cal A}
\subset {\Bbb C}[x_1,x_2,x_3,x_4]$ is minimally
generated by one quadric and 
$r-1$ binomials of degree $r-1$.
In this example $X_{\Cal A}$ is not Cohen-Macaulay. \qed
\endexample

A well-known result due to M.~Hochster states that
normal affine toric varieties are Cohen-Macaulay.
A consequence of this fact is the following degree bound.

\proclaim{Theorem 2.7}
If $Y_{\Cal A}$ is a projectively normal
$r$-dimensional toric variety then
the homogeneous toric ideal $I_{\Cal A}$ is generated by
binomials of degree at most $r$.
\endproclaim

This degree bound is sharp since the toric hypersurface
$\,\,x_0^r = x_1 x_2 \cdots x_r \,\,$ is projectively normal.
It is unknown whether Theorem 2.7 extends to Gr\"obner bases.
In the following conjecture we do not allow any linear changes of coordinates.

\proclaim{Conjecture 2.8}
If $Y_{\Cal A}$ is projectively normal
$r$-dimensional toric variety, then $I_{\Cal A}$ has
a Gr\"obner basis consisting of
binomials of degree at most $r$.
\endproclaim

To appreciate the distinction between generators and 
Gr\"obner bases consider the following example:
If $d=3,n=8$ and ${\Cal A} = \{
(3,0,0), (0,3,0), (0,0,3), (2,1,0),$ \break $
(1,2,0), (2,0,1), (1,0,2), (0,2,1), (0,1,2) \}$,
then $I_{\Cal A}$ is generated by quadrics
but has no quadratic Gr\"obner basis. It is unknown 
whether ${\Bbb C}[{\Cal A}]$ is a Koszul algebra.

Another problem is to better understand the effect 
of smoothness on the degrees of the defining equations.

\proclaim{Conjecture 2.9}
If $Y_{\Cal A}$ is a smooth and projectively normal
toric variety, then the toric ideal $I_{\Cal A}$ is
generated by quadratic binomials.
\endproclaim

It may even be conjectured that in this case
$I_{\Cal A}$ possesses a quadratic Gr\"obner basis.
It was briefly believed in January '95 that Conjecture 2.9 had
been proven, but that proof was withdrawn.
The answer is affirmative for scrolls by 
a result of Ewald and Schmeinck \cite{EwS}.
Conjecture 2.9 is also known to be true for
toric surfaces. This is a consequence of the following theorem
due to Bruns, Gubeladze and Trung \cite{BGT}

\proclaim{Theorem 2.10}
Let $Y_{\Cal A}$ be a projectively normal toric surface
and suppose that the polygon $conv({\Cal A})$ has at
least four lattice points on its boundary.
Then $I_{\Cal A}$ possesses a quadratic lexicographic Gr\"obner basis.
\endproclaim

There is an important sufficient condition for $Y_{\Cal A}$ to 
be projectively normal. (It is not necessary; see Example 13.17
in \cite{Stu}).
A {\it triangulation} $\Delta$ of the configuration 
${\Cal A}$ is a triangulation
of the polytope $conv({\Cal A})$ whose vertices lie in ${\Cal A}$.
We call the triangulation $\Delta$ {\it unimodular} if
each simplex in $\Delta$ is a primitive lattice simplex.

\proclaim{Theorem 2.11}
There exists a unimodular triangulation of ${\Cal A}$
if and only if some initial monomial ideal of $I_{\Cal A}$ is radical.
In this case $Y_{\Cal A}$ is projectively normal.
\endproclaim

The following theorem was established by Knudsen and Mumford in
the early days of toric geometry. It is a key ingredient in Mumford's 
proof of the semi-stable reduction theorem.
See \cite{KKMS} for details and the proof of Theorem 2.12.

\proclaim{Theorem 2.12}
Let $Q$ be any lattice polytope in ${\Bbb Z}^d$.
There exists an integer $m \gg 0$ such that the configuration
$\,m Q \,\cap \, {\Bbb Z}^d$ possesses a unimodular triangulation.
\endproclaim

It would be interesting to find an effective version of this theorem.

\proclaim{Problem 2.13}
Does there exist a bound $M(d)$ such that,
for every $m \geq M(d)$ and  every
lattice polytope $Q$ in ${\Bbb Z}^d$,
the set $m Q \,\cap \, {\Bbb Z}^d $ has
a unimodular triangulation~?
\endproclaim

In our discussion so far $Y_{\Cal A}$ was given
as an explicit subvariety of some projective space
${\Bbb P}^{n-1}$. What can be said about $Y_{\Cal A}$ as an  abstract 
scheme, independently of any choice of a very ample line bundle ?
This is where {\it polyhedral fans} enter the picture.
Let $Q$ be any polytope in a real vector space $V$.
For a face $F$ of $Q$ we consider the set of
linear functionals on $Q$ which attain their maximum at $F$.
This is a convex polyhedral cone in the dual space $V^*$.
The collection of these cones, as $F$ runs over all faces of
$Q$, is a polyhedral fan. It is called the {\it normal fan} of $Q$.

\proclaim{Proposition 2.14}
Two projective toric varieties $Y_{\Cal A}$
and $Y_{\Cal B}$ have isomorphic normalizations
if and only if their polytopes $conv({\Cal A})$
and $conv({\Cal B})$ have the same normal fan.
\endproclaim

Here ``isomorphic'' refers to an isomorphism  of equivariant torus embeddings.
If $Y_{\Cal A}$ is normal, then the normal fan of $conv({\Cal A})$
retains just enough information to remember $Y_{\Cal A}$
as an abstract torus embedding. But it forgets the specific 
line bundle which was used to map $Y_{\Cal A}$ into
${\Bbb P}^{n-1}$. In the synthetic approach to toric 
varieties, as presented in the books \cite{Ful} and \cite{Oda},
the normal fan comes before the polytope. Starting with any
complete fan, one constructs an abstract complete normal 
toric variety by gluing affine pieces 
as in Proposition 2.4.
See also \cite{Cox}.

\example{Example 2.15}
Fix positive integers $i<j<k$ and consider the configuration
$$ {\Cal A} \quad := \quad conv(\{ (i,j,k), (i,k,j), (j,i,k),
(j,k,i), (k,i,j), (k,j,i)\}) \,\cap \,{\Bbb Z}^3 .$$
The normal fan of the hexagon $conv({\Cal A})$ 
is independent of the choice of $i,j,k$.
All toric surfaces $Y_{\Cal A}$ arising
for different choices of $i,j,k$ are normal
and isomorphic to one another. The abstract scheme
$Y_{\Cal A}$ equals ${\Bbb P}^2$ blown up at three points.
All ideals $I_{\Cal A}$ possess 
quadratic lexicographic Gr\"obner bases, by Theorem 2.10.
\endexample

In many applications of toric geometry one encounters
toric varieties $Y_{\Cal A}$ which are not normal.
Or sometimes they are normal but this is difficult
to verify. We describe two instances of the latter
kind arising from representation theory.

Let $G$ be a connected semi-simple
algebraic group over ${\Bbb C}$. Fix
a maximal torus $({\Bbb C}^*)^d$ in $G$, let $P$ be a 
parabolic subgroup containing $({\Bbb C}^*)^d$
and consider the flag variety $G/P$.
The following result is due to R.~Dabrowski \cite{Dab}.

\proclaim{Theorem 2.16} 
The closure in $G/P$ of a generic $({\Bbb C}^*)^d$-orbit is normal.
\endproclaim

It is open whether the generic torus orbit closures are 
{\it projectively normal} for all very ample line bundles on  $G/P$.
It is also open whether the closures of 
all (non-generic) $({\Bbb C}^*)^d$-orbits are normal,
or even projectively normal.
This latter  conjecture is known to be true
in the case  when $G/P$ is the classical Grassmannian 
$\,Gr_r ({\Bbb C}^d)\,$ in its Pl\"ucker embedding.
In this case ${\Cal A}$ is a subset of the {\it hypersimplex}
$$ \bigl\{ (i_1,\ldots,i_d) \in \{0,1\}^d \,:
\,i_1 + \cdots + i_d = r \bigl\} $$
such that ${\Cal A}$ consists of the incidence vectors of 
the bases of a realizable 
{\it matroid}. For instance, in Example 1.2 (d)
that matroid is the uniform rank $2$ matroid on $4$ elements.
We refer to \cite{GGMS} for an introduction to matroids from
an algebro-geometric point of view.
The following result is due to Neil White \cite{Wh1}.

\proclaim{Theorem 2.17} 
Let ${\Cal A}$ be the set of incidence vectors of the
bases of a matroid. Then  $Y_{\Cal A}$ is projectively normal.
\endproclaim

We close with a reformulation of a combinatorial conjecture in \cite{Wh2}.

\proclaim{Conjecture 2.18} 
Let ${\Cal A}$ be the set of incidence vectors of the bases of a matroid. 
Then the homogeneous toric ideal $I_{\Cal A}$ is generated by quadratic
binomials.

\endproclaim

\vskip 1cm

\head 3. Binomial Zoo \endhead

Our objective is to understand the minimal generators and
Gr\"obner bases of the toric ideal $I_{\Cal A}$.
To this end we introduce the following three definitions.
A binomial ${\bold x}^{{\bold u}_+} - {\bold x}^{{\bold u}_-} $
in $I_{\Cal A}$ is called a {\it circuit} if its support
(i.e.~the set of variables appearing in that binomial) is 
minimal with respect to inclusion.
We write ${\Cal C}_{\Cal A}$ for the set of all circuits in $I_{\Cal A}$.
Geometrically speaking, we consider all images of
$X_{\Cal A}$ under projection into
coordinate subspaces of ${\Bbb P}^{n-1}$;
such an image is called a circuit of $X_{\Cal A}$ if
it has codimension $1$ in its coordinate subspace.
We define the {\it universal Gr\"obner basis} 
$\,{\Cal U}_{\Cal A}$ to be the union of all reduced 
Gr\"obner bases of $I_{\Cal A}$. We say that a binomial
${\bold x}^{{\bold u}_+} - {\bold x}^{{\bold u}_-} $ in $I_{\Cal A}$
lies in the {\it Graver basis} $Gr_{\Cal A}$ if there exists no other binomial
${\bold x}^{{\bold v}_+} - {\bold x}^{{\bold v}_-} \in I_{\Cal A}$
such that ${\bold x}^{{\bold v}_+}$ divides ${\bold x}^{{\bold u}_+}$
and ${\bold x}^{{\bold v}_-}$ divides ${\bold x}^{{\bold u}_-}$.

\proclaim{Proposition 3.1}
For any toric ideal $I_{\Cal A}$ we have the inclusions
$$ {\Cal C}_{\Cal A} \,\subseteq \, 
{\Cal U}_{\Cal A} \,\subseteq \,  Gr_{\Cal A} . \tag 3.1 $$
Each of the four combinations of strict or non-strict inclusions is possible.
\endproclaim

\example{Examples 3.2}
\roster
\item
For the twisted cubic curve in
Example 1.2 (a) we have 
$$ \align
{\Cal C}_{\Cal A} \,& =\,  \{
x_1 x_3 - x_2^2, \,
x_2 x_4 - x_3^2, \,
x_1^2 x_4 - x_2^3, \,
x_1 x_4^2 - x_3^3 \}  \quad \hbox{and} \quad \cr
{\Cal U}_{\Cal A} \,&=\,\, Gr_{\Cal A} 
\,\, = \,\, {\Cal C}_{\Cal A} \,\cup \,\{ x_1 x_4 - x_2 x_3 \}.
\endalign $$
In this example the circuits do not generate
the ideal $I_{\Cal A}$ although they do define
the twisted cubic as a subscheme of ${\Bbb P}^3$.
\item
For the Veronese surface in ${\Bbb P}^5$,
the set of circuits ${\Cal C}_{\Cal A}$ equals
the universal Gr\"obner basis ${\Cal U}_{\Cal A}$, but
the Graver basis $Gr_{\Cal A}$ properly contains ${\Cal U}_{\Cal A}$.
\endroster
\endexample

In general, the circuits of a projective toric variety 
$Y_{\Cal A}$ do not define $Y_{\Cal A}$
scheme-theoretically. But they do set-theoretically.
We state this result in the affine case.

\proclaim{Theorem 3.3}
The toric variety $X_{\Cal A}$ is cut out by its circuits:
$Rad \langle {\Cal C}_{\Cal A} \rangle = I_{\Cal A}$.
The embedded components of $\langle {\Cal C}_{\Cal A} \rangle$
are supported on torus orbits of $X_{\Cal A}$,
i.e.~every associated prime of $\langle {\Cal C}_{\Cal A} \rangle$
has the form $I_{F \cap {\Cal A}}$ for some face $F$ 
of the cone $pos({\Cal A})$.
\endproclaim 

Theorem 3.3 was proved jointly with David Eisenbud in \cite{EiSt}.
The main result in \cite{EiSt}
concerns the decomposition of  arbitrary binomial schemes.
It underlines the importance of toric ideals
as combinatorial building blocks in algebraic geometry:

\proclaim{Theorem 3.4}
Let $I$ be any binomial ideal in ${\Bbb C}[x_1,\ldots,x_n]$.
Then $I$ possesses a primary decomposition into primary
binomial ideals. In particular, every associated prime
of $I$ has the form $\,I_{\Cal A} + 
\langle x_{i_1},\ldots,x_{i_r} \rangle \,$ for some
configuration ${\Cal A}$.
\endproclaim

There is an important class of toric varieties for which
both inclusions in (3.1) are equalities. They are called
{\it unimodular} and defined by the following theorem.

\proclaim{Theorem 3.5}
For a configuration ${\Cal A}$ the following conditions are equivalent:
\roster
\item
Every triangulation of ${\Cal A}$ is unimodular.
\item
For every subset ${\Cal B}$ of ${\Cal A}$, the quotient
group ${\Bbb Z}{\Cal A}/{\Bbb Z}{\Cal B}$ is free abelian.
\item
Every initial monomial ideal of $I_{\Cal A}$
is a radical ideal.
\item
In every circuit of $I_{\Cal A}$ both monomials are square-free.
\endroster
\endproclaim

If these four equivalent conditions hold, then we call ${\Cal A}$
{\it unimodular}. 

\proclaim{Proposition 3.6}  
If ${\Cal A}$ is unimodular then $X_{\Cal A}$ is normal.
\endproclaim

\example{Example 3.7}  
The term ``circuits'' originates from the following class of
unimodular toric varieties. Let $G$ be a finite directed graph.
We label the vertices of $G$ by $1,2,\ldots,d$ and we encode
the edges of $G$ as differences of unit vectors in ${\Bbb Z}^d$:
$$ {\Cal A}_G \quad = \quad \bigl\{\, {\bold e}_i - {\bold e}_j 
\in {\Bbb Z}^d \,:\,
(i,j) \, \hbox{is a directed edge of $G$} \,\bigr\}. $$
The unimodularity of ${\Cal A}_G$ is a basic result in matroid theory.
The circuits in the toric ideal 
$I_{{\Cal A}_G}$ correspond to the directed circuits in $G$. The
following examples of circuits of length five illustrate this correspondence:
$$ \matrix
 \hbox{Directed circuit in $G$} & &
\leftrightarrow & & \hbox{Circuit in $I_{{\Cal A}_G}$} \\
1 \rightarrow
2 \rightarrow
3 \rightarrow
4 \rightarrow
5 \rightarrow 1 & &
\leftrightarrow  & & x_{12} x_{23} x_{34} x_{45} x_{51} - 1 \\
1 \rightarrow
2 \rightarrow
3 \leftarrow
4 \rightarrow
5 \rightarrow 1 & &
\leftrightarrow  & & x_{12} x_{23} x_{45} x_{51} -  x_{43}  \\
1 \rightarrow
2 \rightarrow
3 \leftarrow
4 \rightarrow
5 \leftarrow 1 & &
\leftrightarrow  & & x_{12} x_{23} x_{45} - x_{15} x_{43}  \\
\endmatrix $$
If $G$ is the complete directed bipartite graph $K_{r,s}$ then 
$I_{{\Cal A}_G}$ is the ideal of $2 \times 2$-minors of an
$r \times s$-matrix of indeterminates
and $Y_{{\Cal A}_G}$ is the Segre variety $ {\Bbb P}^{r-1} 
\times {\Bbb P}^{s-1} $. The degrees of the circuits in
$I_{{\Cal A}_G}$ range from $2$ to $max  \{r,s \}$.

If $G$ is a complete graph, then ${\Cal A}_G$ equals
the {\it root system} of type $A_{n-1}$. The toric
variety $X_{{\Cal A}_G}$ is the closure of a generic
$({\Bbb C}^*)^n$-orbit in the {\it adjoint representation}
of $GL_n({\Bbb C})$. Passing to subgraphs corresponds to
passing to non-generic orbit closures in the adjoint representation. 
Thus the affine toric varieties defined by directed graphs are precisely 
the closures of $({\Bbb C}^*)^n$-orbit in the adjoint representation
of $GL_n({\Bbb C})$. These toric varieties are all normal
by Proposition 3.6. \qed
\endexample

A configuration $X_{\Cal A}$ is called {\it hereditarily normal}
if all closures of torus orbits on  $X_{\Cal A}$ are normal.
This is equivalent to the requirement that
$X_{\Cal B}$ is normal for every subset ${\Cal B}$ of ${\Cal A}$.
Unimodular configurations are hereditarily normal, but not conversely.
Normality of torus orbits for abitrary semisimple algebraic groups
was studied by J.~Morand \cite{Mor}. She showed that the
root systems of types $A_n,\; B_2,\; C_2,\;D_4$ are hereditarily normal
and that all other root systems are not hereditarily normal.

The following criterion underlines
the geometric significance of the circuits.

\proclaim{Theorem 3.8}
An affine toric variety $X_{\Cal A}$ is hereditarily normal
if and only if every circuit in $I_{\Cal A}$ has at least one
square-free monomial.
\endproclaim

The Graver basis $Gr_{\Cal A}$ has the following
geometric interpretation. Consider the closure
of the affine variety $X_{\Cal A} \subset {\Bbb C}^n$
in the $n$-fold product of projective lines 
$\,{\Bbb P}^1 \times {\Bbb P}^1 \times \cdots \times {\Bbb P}^1  $.
It is defined by a toric ideal $ \, I_{\Lambda({\Cal A})} \,$ in
${\Bbb C}[x_1,y_1, x_2,y_2, \cdots, x_n,y_n]$.
which is homogeneous with respect to each pair
of variables $(x_i,y_i)$.
Here $\Lambda({\Cal A})$ is a certain configuration of $2n$ vectors
in ${\Bbb Z}^{n+d}$, which is called the {\it Lawrence lifting} of
${\Cal A}$. A key property of such configurations $\Lambda({\Cal A})$
is the following.

\proclaim{Theorem 3.9}
The defining ideal $I_{\Lambda({\Cal A})}$ of a toric subvariety of
$\,{\Bbb P}^1  \times \cdots \times {\Bbb P}^1  \,$ is minimally 
generated by its Graver basis $Gr_{\Lambda({\Cal A})}$.
\endproclaim

The original ideal $I_{\Cal A}$ is recovered from $I_{\Lambda({\Cal A})}$ 
by dehomogenizing, that is, by replacing the
variables  $\,y_1,\ldots,y_n \, $ by $ 1$.
This induces a bijection between the Graver basis of
$\Lambda({\Cal A}) $ and the Graver basis of ${\Cal A}$.
Under this bijection we have the following 
geometric interpretation of the Graver basis.

\proclaim{Corollary 3.10}
The Graver basis of an affine toric variety
consists of the minimal generators of the ideal
of its closure in $\,{\Bbb P}^1 \times  \cdots \times {\Bbb P}^1 $.
\endproclaim

For instance, if ${\Cal A}$ is the configuration in Example 1.2 (a), then
$$ \align
 I_{\Lambda({\Cal A})} \quad = \quad
 \langle \,
   & x_1 x_3 y_2^2 - x_2^2 y_1 y_3 , \,\,
   x_1 x_4 y_2 y_3 - x_2 x_3 y_1 y_4 , \,\,  
   x_2 x_4 y_3^2 - x_3^2 y_2 y_4 , \,   \\
   &\,\, x_1^2 x_4 y_2^3- x_2^3 y_1^2 y_4 ,\,\,
   x_1 x_4^2 y_3^3 - x_3^3 y_1 y_4^2 \,\rangle. \\ 
 \endalign $$
This ideal defines a toric surface in
${\Bbb P}^1 \times {\Bbb P}^1 \times {\Bbb P}^1 \times {\Bbb P}^1$.
It is the closure of the affine cone over the twisted cubic curve.
We get the Graver basis for $I_{\Cal A}$,
the ideal of the twisted cubic, by setting
$y_1,y_2,y_3,y_4$ to $1$ in these five binomials.

\vskip 1cm

\head 4. Degree Bounds \endhead

The following problem is of great interest in
computational algebraic geometry.

\proclaim{Conjecture 4.1}
Let $Y$ be an irreducible projective variety of degree $d$. Then the 
homogeneous prime ideal of $Y$ is generated by forms of degree at most $d$.
\endproclaim

Many authors have concentrated on proving
stronger inequalities involving the {\it
(Castelnuovo-Mumford) regularity} $reg(Y)$ of $Y$.
Let $I_Y$ denote the ideal of $Y$ in ${\Bbb C}[x_1,\ldots,x_n]$.
Then $reg(Y)$ is the maximum of the numbers 
$deg( \sigma)-i$ where $\sigma$ is a minimal $i$-th syzygy of $I_Y$. 
Here $i=0$ is for ideal generators,
$i=1$ for syzygies among them, etc...
Following \cite{BS},
two equivalent definitions of {\it regularity} are:
\roster
\item $reg(Y)$ is the maximum degree of an element in 
the reduced  Gr\"obner basis of $I_{\phi(Y)}$ with respect
to the reverse lexicographic term order, where 
$I_{\phi(Y)}$ is the image of $I_Y$ 
under a generic linear automorphism $\phi \in GL_n({\Bbb C}) $;

\item $reg(Y)$ is the smallest integer such that
$\,H^i( {\Bbb P}^{n-1}, {\Cal I}_Y (s)) \, = \,0 \,$ for all
$i \geq 0 $ and $s \geq reg(Y) - i$, where
${\Cal I}_Y$ is the ideal sheaf of $Y$.
\endroster

\vskip .2cm

\noindent
The following conjecture due to Eisenbud and Goto \cite{EG} implies
Conjecture 4.1.

\proclaim{Conjecture 4.2}
Let $Y$ be an irreducible projective variety,
not contained in any hyperplane. Then
$\,\, reg(Y) \,\, \leq \,\, degree(Y) - codim(Y) + 1 $.
\endproclaim

The Eisenbud-Goto inequality holds when $Y$ is 
arithmetically Cohen-Macaulay. Hence, by Hochster's Theorem, 
it holds for projectively normal toric varieties.
Conjecture 4.2 was proved for curves by Gruson-Lazarsfeld-Peskine 
\cite{GLP}, and for irreducible smooth surfaces and $3$-folds by 
Lazarsfeld \cite{Laz} and Ran \cite{R}.

Conjectures 4.1 and 4.2 are widely open in general, even for 
toric varieties $Y_{\Cal A}$. Recall from Theorem 1.4 that the 
degree of $Y_{\Cal A}$ equals the volume of the polytope
$\,conv({\Cal A})$. Also the regularity of $Y_{\Cal A}$
can be expressed combinatorially, using the
simplicial representation of {\it Koszul homology}
(see e.g.~Theorem 12.12 in \cite{Stu}). 
A class of toric varieties for which the Eisenbud-Goto inequality
is valid was identified by Irena Peeva and the author in \cite{PS}.
For toric varieties $Y_{\Cal A}$ in codimension 2 we 
explicitly construct the minimal free resolution of $I_{\Cal A}$.
Our construction implies:

\proclaim{Theorem 4.3}
Let $Y_{\Cal A}$ be a projective toric variety of codimension $2$ in 
${\Bbb P}^{n-1}$, not contained in any hyperplane.
Then $\,reg(Y_{\Cal A}) \,\leq \, degree(Y_{\Cal A}) - 1 $.
\endproclaim 

\example{Example 4.4}
The following family of toric surfaces in ${\Bbb P}^4$ shows
that the inequality in Theorem 4.3 is tight. For any
integer $d \geq 3$ let
$$ {\Cal A} \quad = \quad
\pmatrix
 1 &  1 &  1 &  1 &  1 \\
 0 &  1 &  1 &  0 &  0 \\
 0 &  0 &  1 &  1 &  d \\
\endpmatrix .$$
The surface $Y_{\Cal A}$ is not arithmetically
Cohen-Macaulay,
its toric ideal $I_{\Cal A}$ is generated by one quadric
and $d$ forms of degree $d$, and
$\,degree(Y_{\Cal A}) = \hbox{\it Vol}( conv({\Cal A}) ) = d+1 $.
\endexample

Given an arbitrary toric variety in higher codimension, the following
is the currently best general bound for its regularity in terms of its degree.

\proclaim{Theorem 4.5}
Let $Y_{\Cal A}$ be a projective toric variety in ${\Bbb P}^{n-1}$. Then 
$$ reg(Y_{\Cal A}) \quad \leq \quad n  \cdot
degree(Y_{\Cal A}) \cdot codim( Y_{\Cal A}) . \tag 4.1  $$
\endproclaim 

We shall outline the proof of Theorem 4.5. First, we make use of the fact 
that regularity is upper semi-continuous with respect
to flat families. This implies 
$$ reg(\,I_{\Cal A} \,) \quad \leq \quad reg \bigl(
in_\prec(I_{\Cal A}) \bigr) 
\tag 4.2 $$
for any term order $\prec$ on ${\Bbb C}[x_1,\ldots,x_n]$.
The {\it Taylor resolution} (see page 439 in \cite{Eis})
for the monomial ideal
$\,in_\prec(I_{\Cal A}) \,$ implies the inequality
$$ reg( \,in_\prec(I_{\Cal A}) \,) \quad \leq \quad
n \cdot maxgen(\, in_\prec(I_{\Cal A}) \,), \tag 4.3 $$
where  $\,maxgen(\,\cdot \,)\,$ denotes the maximal 
degree of any minimal generator. In view of the
inclusion ${\Cal U}_{\Cal A} \subseteq Gr_{\Cal A}$
in Proposition 3.1, it suffices to prove the following:

\proclaim{Lemma 4.6}
Let $Y_{\Cal A}$ be a projective toric variety. Every binomial in the 
Graver basis $Gr_{\Cal A}$ of the homogeneous ideal $I_{\Cal A}$ has 
degree at most $\, degree(Y_{\Cal A}) \cdot codim( Y_{\Cal A}) $.
\endproclaim

For a set $S$ of polynomials let $\,maxdeg( S ) \,$ denote the maximum 
degree of any element in $ S $. The proof of Lemma 4.6
is derived from the next two inequalities:
$$ maxdeg({\Cal C}_{\Cal A}) \quad \leq \quad degree(Y_{\Cal A}) \tag 4.4 $$
This holds because the hypersurface defined by any circuit
is a projection of $Y_{\Cal A}$. The second inequality is
$$ maxdeg(Gr_{\Cal A}) \quad \leq \quad
codim(Y_{\Cal A}) \cdot maxdeg({\Cal C}_{\Cal A}) .\tag 4.5 $$
This follows from a standard inequality for Hilbert bases of 
integer monoids. See Chapter 4 in \cite{Stu} for details.

While the inquality (4.4) is tight, there
seems to be some room for improvement left in (4.5).
In my lecture at Santa Cruz I asked whether
even the equality $\, maxdeg(Gr_{\Cal A}) \,=\, 
maxdeg({\Cal C}_{\Cal A}) \,$ might be true.
This would have implied Conjecture 4.1 for all toric varieties.
Unfortunately this equality was too good to be true.
Serkan Hosten and Rekha Thomas found the following counterexample.

\example{Example 4.7}
We shall present a projective toric variety in ${\Bbb P}^9$ which satisfies
$$ 16 \,=\,  maxdeg(Gr_{\Cal A}) \, = \,
maxgen(I_{\Cal A}) \,\,\,> \,\,\,maxdeg({\Cal C}_{\Cal A}) = 15.
\tag 4.6 $$
Our starting point is the affine toric surface 
$\,X_{\Cal B} \subset {\Bbb C}^5 \,$ defined by 
$$ {\Cal B} \quad = \quad 
\pmatrix
  1 &  3 & 4 &  \phantom{-} 6 & 0 \\
  0 &  0 & 0 & - 5 & 1 \\
\endpmatrix.
$$
We then construct their Lawrence lifting 
$\,{\Cal A} := \Lambda({\Cal B})$. In other words,
$I_{\Cal A}$ is the toric ideal defining the closure of $\,X_{\Cal B}\,$ in 
${\Bbb P}^1 \times {\Bbb P}^1 \times {\Bbb P}^1 \times
{\Bbb P}^1 \times {\Bbb P}^1 $. This ideal equals
$$
\align & I_{\Cal A}  \, = \, \langle \,
 x_2 y_1^3 - x_1^3 y_2 , \,
x_3 y_1^4 - x_1^4 y_3 , \,
x_3^3 y_2^4 - x_2^4 y_3^3 ,\,
x_4 x_5^5 y_2^2 - x_2^2 y_4 y_5^5 ,\,
 x_4 x_5^5 y_1^6 - x_1^6 y_4 y_5^5 ,\\
& \,\,
\underline{x_4^2 x_5^{10} y_3^3 - x_3^3 y_4^2 y_5^{10}  },\,\,
x_4 x_5^5 y_1^2 y_3 - x_1^2 x_3 y_4 y_5^5,\,
x_4 x_5^5 y_1^3 y_2 - x_1^3 x_2 y_4 y_5^5,\,
x_3 y_1 y_2 - x_1 x_2 y_3, \\
&
  x_1 x_4 x_5^5 y_2 y_3 \! - \! x_2 x_3 y_1 y_4 y_5^5,
  x_2^2 x_4 x_5^5 y_3^3 - x_3^3 y_2^2 y_4 y_5^5 ,
  x_1 x_3^2 y_2^3 - x_2^3 y_1 y_3^2 ,
  x_1^2 x_3 y_2^2 - x_2^2 y_1^2 y_3 ,\\
& \,\,
   x_2 x_4 x_5^5 y_1 y_3^2 - x_1 x_3^2 y_2 y_4 y_5^5 ,\,
  x_1^2 x_4 x_5^5 y_3^2 - x_3^2 y_1^2 y_4 y_5^5 ,\,\,
  x_4^2 x_5^{10} y_1 y_2 y_3^2 - x_1 x_2 x_3^2 y_4^2 y_5^{10} \,
\rangle. \\
\endalign 
$$
By Theorem 3.9, these $16$ minimal generators of $I_{\Cal A}$ 
coincide with the Graver basis $Gr_{\Cal A}$.
The first six binomials in this list are the circuits.
Thus (4.6) holds for this example. We remark that
the projective variety $Y_{\Cal A} \subset {\Bbb P}^9$
defined by $I_{\Cal A}$ satisfies
$codim(Y_{\Cal A}) = 3$,
$degree(Y_{\Cal A}) = 54 $, and
$ reg(Y_{\Cal A}) = 17 $. \qed
\endexample

\vskip .1cm

In light of this counterexample, I wish to propose the
following improvement of (4.5).
Consider any circuit $\,C \in {\Cal C}_{\Cal A}\,$ and  regard
its support $supp(C)$ as a subset of ${\Cal A}$.
The lattice  $\,{\Bbb Z}(supp(C)) \,$
has finite index in the (possibly bigger) lattice
$\,{\Bbb R}(supp(C)) \cap {\Bbb Z}{\Cal A}$.
We call this index the {\it index} of the circuit $C$.
We define the {\it true degree} of the circuit $C$
to be the product $\,degree(C) \cdot index(C) $.
It can be shown that the true degree of any circuit
is bounded above by $degree(Y_{\Cal A})$.
The following conjecture would thus imply
Conjecture 4.1 for projective toric varieties.

\proclaim{Conjecture 4.8}
The degree of any element in the Graver basis
$Gr_{\Cal A}$ of a toric ideal $I_{\Cal A}$ is bounded above  
by the maximal true degree of any circuit $C \in {\Cal C}_{\Cal A}$.
\endproclaim

Note that Conjecture 4.8 is consistent with Example 4.7
because the underlined circuit
$\, x_4^2 x_5^{10} y_3^3 - x_3^3 y_4^2 y_5^{10}  \,$
has index two and hence its true degree is $30$.

\vfill \eject

\enddocument